\documentclass{amsart}[12 pt]

\usepackage{latexsym}
\usepackage[utf8]{inputenc}
\usepackage{amsmath,amssymb,dsfont,amsfonts}
\usepackage{mathrsfs}
\usepackage{verbatim}
\usepackage{graphicx}
\usepackage{graphics}
\usepackage{geometry}
\usepackage{booktabs}
\usepackage{enumitem}
\usepackage{subfigure}
\usepackage[ruled, vlined, linesnumbered]{algorithm2e}
\usepackage[final]{pdfpages}
\usepackage{color}

\usepackage[colorlinks]{hyperref}

\newcommand{\Prob}{\mathbb P}
\newcommand{\R}{\mathbb{R}}

\newcommand{\E}{\mathbb{E}}

\newcommand{\eop }{ \hfill $\Box$ }

\newtheorem{theorem}{Theorem}[section]
\newtheorem{definition}{Definition}[section]

\newtheorem{lemma}{Lemma}[section]

\newtheorem{assum}{Assumption}

\DeclareMathOperator*{\esssup}{ess\,sup} 
\DeclareMathOperator*{\argmin}{arg\,min} 

\numberwithin{equation}{section}

\allowdisplaybreaks

\begin{document}

%
%

\markboth{}{}

\title[Stochastic near-optimal control: additive, multiplicative, non-Markovian and applications] {Stochastic near-optimal control: additive, multiplicative, non-Markovian and applications}

\author{Lourival Lima}
\address{Mathematics Department \\ Rua Sérgio Buarque de Holanda, 651\\ UNICAMP - State University of Campinas\\ 13083-859 Campinas, Brazil \\Research supported by CAPES 88882.329064/2019-01.}\email{l190685@dac.unicamp.br}

\author{Paulo Ruffino}
\address{Mathematics Department \\ Rua Sérgio Buarque de Holanda, 651\\ UNICAMP - State University of Campinas\\ 13083-859 Campinas, Brazil \\Research partially supported by CNPq 305212/2019-2, FAPESP 2020/04426-6 and 2015/50122-0.}\email{ruffino@unicamp.br}

\author{Francys Souza}
\address{Mathematics Department \\ Rua Sérgio Buarque de Holanda, 651\\ UNICAMP - State University of Campinas\\ 13083-859 Campinas, Brazil \\Research supported by FAPESP 2017/23003-6.}\email{francys@estatcamp.com.br}

\date{November 2020}

\begin{abstract}
In this survey we present the near-optimal stochastic control problem according to some recent tools in the literature. In particular, we focus on the approach of a discretization of the noise values instead of the canonical time-discretization. This is the so called {\it skeleton} structure. This allows to obtain an $\epsilon$-optimal control in non-Markovian systems (the main Theorem). A simple example illustrates the technique. The importance of the approach is emphasised in a final section on open problems related to more geometrical framework and discontinuous noise.
\end{abstract}

\maketitle

\section{Introduction}

 Control theory plays a major role in most applications of differential equation in any physical system. This theory is crucial when one can control -- either constant in time or time-dependent, one or more parameters of a system, say, with conditions like: temperature, pressure, concentration of substances, investments, humidity, position and velocity of autonomous vehicles, satellites, electromagnetic parameters, action with vaccination in a population etc, just to mention few of them. Around the last few decades of the 20th century, emboldened by the well development of deterministic control theory and the constant improving of stochastic analysis and stochastic dynamics, the theory of stochastic control started to develop rapidly thanks also to a countless number of relevant application. As for the deterministic control theory, among hundreds of excellent introductory literature, we mention e.g. Colonius and Kliemann \cite{Col-Kliemann}, Bullo and Lewis \cite{Bullo-Lewis}, Bacciotti \cite{Bacciotti} and references therein; yet, for stochastic analysis, dynamics and control, among a list  of excellent introductory texts, see e.g. Arnold \cite{Arnold}, Oksendal \cite{Oksendal}, Protter \cite{Protter} and references therein.

Control theory in stochastic systems is fascinating in the sense that although the outcome is random and unpredictable, nevertheless in many cases, its law as a random variable can be controlled. It means that many useful properties and tools can be applied in order to optimise the chance that the outcomes are favourable. The purpose of this review is to show some of these tools and in particular to show an application which describes an algorithm to obtain an $\epsilon$-sub-optimal control. The open problems in the final section show how challenging and significant this topic is.

This review is organised as follows: after presenting the general set-up in the Introduction, in Section 2 we introduce the main framework and discretization structure. Section 3 presents recent tools and results concerning near optimality (Theorem 1). A toy model illustrates the algorithm of the main theorem. Section 4 summarises some open relevant questions in a wider scenario of near optimal stochastic control: we discuss the importance of the multiplicative equation and the case of discontinuous noise. Namely: The first open problem regards stochastic control systems in a Lie group, intended to include matrix product; the second problem aims L\'evy noise.

\subsection{The set-up and the model}

We consider a filtered probability space $(\Omega, \mathcal{F}, \mathbb{F}, \Prob)$, where the $\sigma$-algebra filtration $\mathbb{F}:=\mathcal{F}_t$ is generated by a stochastic process $Z_t$, $t\in [0,T]$ (the noise), which is a functional of the Brownian motion in the real line, including non-Markovian and non-semimartingale processes controlled by mutually singular measures. We shall denote by $X^u (t)$ the trajectory (the outcome) in the state space corresponding to a certain control function $u$. Our model includes functional dependence on the past trajectory (see more related to pseudo-Markovian equation in \cite{Claisse-Talay-Tan}). That is: at any time $t\in [0,T]$, the whole past trajectory of the system in $[0,t]$, denoted by the sub-index $X
^u_ t$ is considered in the equation. Hence, we have $X^u_ t:\Omega \times [0,t] \rightarrow \R^n  $ in such a way that $X^u_t(s)= X^u_{t'}(s)$ for all $s\leq t\leq t'$, in particular $X^u_t(t)= X^u(t)$. We assume that $ X
^u (t)$  is modelled by the following functional stochastic differential equation:
\begin{equation}   \label{eq: model}
dX^u (t) = \alpha(t, X^u_ t, u(t))\, dt + \sigma(t, X^u_t, u(t))\, dZ(t),
\end{equation}
where $\alpha$ and $\sigma$ are vector fields sufficiently regular depending on $t$, on the past trajectories $X^u_t$ and on the control function $u$. More precisely, in our non-antecipative model, the control function $u: \Omega \times [0,T] \rightarrow \R^k$ lives in the space $U^T$ of bounded, $\mathcal{F}_t$-adapted stochastic processes. At each time $t\in [0,T]$, the "best future" control segment will be looked for in the space  $U^T_t; 0\le t \le T$, a suitable family of admissible $\mathbb{F}$-adapted controls defined over $(t,T]$, such that  $U^T:=U_0^T$.

The aim of the controller is to optimise the expectation of a given functional $\xi$, say, the payoff, energy, cost, distance etc. More precisely, denote by $\mathbf{C}([0,T])$ the set of continuous functions from $[0,T]$ to the state space $\R^n$, then $\xi:\mathbf{C}([0,T])\rightarrow \mathbb{R}$ is a Borel ({\it payoff}, say) functional. Our aim is, theoretically,  to find a control function $u(t)$ which achieve the optimal successful result:

\begin{equation*}
\sup_{u\in U^T} \E \left[  \xi (X^u)\right].
\end{equation*}
It is well known in the literature that such an optimal control function may not exist, see e.g. \cite{LOS} or \cite{Zhou}. In view of this fact, we look for an algorithm which provides a sufficiently close approximation in the sense that given any $\epsilon>0$, it is possible to construct a control $u^*$ such that
\begin{equation}\label{cp1}
\E \left[  \xi (X^{u^*})\right] > \sup_{u\in U^T} \E \left[  \xi (X^u)\right] -\epsilon.
\end{equation}

The original main problem of calculating an optimal control is essentially an infinite-dimensional nonlinear optimisation question in which explicit solutions are not available. Our strategy is to reduce, via a convenient type of discretization, into a finite dimensional problem, such that near-optimal (i.e. $\epsilon$-optimal) adapted controls satisfying equation (\ref{cp1}) can be calculated. A proof of convergence of the discretized version to the original continuous problem can also be done. Note that in this model the functional $\xi$ depends only on the state space. This is not quite a big restriction since increasing the dimension of the state space $\R^n$ one can include dependence of $\xi$ also of the control $u(t)$, on time $t$, or on any other relevant parameter of the model.

\section{Classical results and new tools}

In this section, we summarise some of  the classical results in the literature concerning stochastic optimal control and respective numerical methods. In the Markovian case, the traditional approach is based on Hamilton-Jacobi-Bellman (HJB) PDEs and based in discretizations into Markov chains, see e.g. \cite{krylov1, kushner2}. We also remember that the use of PDE techniques is only numerically efficient in low dimensions.

In addition to the Markovian context, the value process cannot be characterised by HJB-type PDE's, and consequently, the solution of the stochastic optimal control problem is much more demanding. For details, see e.g., among many others \cite {nutz2, qiu, tan, saporito, zhang, ren}. The PDE approach in the sense of \cite{dupire, cont} for the problem of non-Markovian control in continuous time, in general, is restricted to the theoretical characterisation of the value process, not allowing to derive a numerical algorithm for the concrete resolution of the problem, that is, finding optimal or approximately optimal strategies. Such an approach is not computationally treatable, especially in situations in which the state of the system is high-dimensional and non-elliptical.

Recently, it has been developed a theory of functional stochastic calculus and control applications in a series of studies \cite{LEAO_OHASHI2013, LEAO_OHASHI2018.1, LEAO_OHASHI2018.2, LEAO_OHASHI2019.1, LOS}. The main advantage of this approach over other methods lies in the numerically viable description of stochastic systems controlled in continuous time, especially in problems of optimal control, possibly non-Markovian and of high dimension. In particular, in \cite{LOS}, the authors develop  a discretization method that allows to calculate near-optimal controls for a given optimal control problem in a range of $ [0, T] $ with $ 0 <T < \infty $, via a discrete-time PDE.

From now on, we describe the relevant aspects of this construction. In few words: the main idea in this new approach is to consider a discretization based on a partition in the values of the noise, instead of a partition in the time line. The control function is then updated towards best performance at each stopping time associated to this partition. The updating of the control is performed using dynamics programming (DP). The advantage of this approach is that it naturally give us a numerical algorithm for simulation (see the numerical example in Section \ref{Numerical Example} below).

For a pair of finite $\mathbb{F}$-stopping times $(M,N)$, we denote the random intervals in the standard way: $ ]]M,N ]]: = \{(\omega,t); M(\omega) < t \le N(\omega)\}$ and $]] M, +\infty[[:=\{(\omega,t); M(\omega)< t < +\infty\}.$ The control function is bounded with values in the so called {\it action space} which is going to be the compact cube

$$A:=\{(x_1, \ldots, x_r)\in \mathbb{R}^r; \max_{1\le i\le r}|x_i|\le a\}$$
for some $0 < a< +\infty$.

In order to set up the basic structure of our control problem, we first need to define the class of admissible control processes: For each pair $(M,N)$ of a.s finite $\mathbb{F}$-stopping times such that $M < N$ a.s, we denote

\begin{equation*}\label{controlset}
U^{N}_{M}:=\{ \mathbb{F} \text{ -predictable processes }~u:~]]M,N]]\rightarrow A;\ \mbox{ such that the limit }  u(M+)\  \text{exists}\}.
\end{equation*}
Note that whenever necessary, we can always extend a given $u\in U^N_M$ by setting $u = 0$ in the complement of the random interval $]]M,N]]$.
Elements in the  family of processes in $U^{N}_{M}$ satisfy the following straightforward properties:

\bigskip

\begin{description}
  \item[1)] \textbf{Restriction}: $u\in U_M^N\Rightarrow u\mid_{]]M,P]]}\in U_{M}^P$ for $M < P \le N$~a.s.

  \bigskip

  \item[2)] \textbf{Concatenation}: If $u\in U_{M}^N$ and $v\in U_N^P$ for $M<N< P~a.s$, then $(u\otimes_N v)(\cdot)\in U_M^P$, where
  \begin{equation*}\label{concatenation}
  (u\otimes_N v)(r):=\left\{
\begin{array}{rl}
u(r); & \hbox{if} \ M < r \le N \\
v(r);& \hbox{if} \ N < r \le P.
\end{array}
\right.
\end{equation*}

\bigskip

  \item[3)] \textbf{Finite composition}: For every $u,v\in U_M^N$ and $K\in \mathcal{F}_M$, we have
$$u\mathds{1}_K + v\mathds{1}_{K^c}\in U^N_M.$$

\end{description}

By  $U_M$ we mean $U_M^{\infty}$. For each $1\le p< \infty$, let $L^p_a(\mathbb{P}\times Leb)$ be the Banach space of all $\mathbb{F}$-adapted $\R^n$-valued processes. I.e.  $Y\in L^p_a(\mathbb{P}\times Leb)$ if and only if it is adapted and

\begin{equation*}\label{Lpspace}
\mathbb{E}\int_0^T\|Y(t)\|^p \ dt < \infty.
\end{equation*}
Consider $\mathbf{B}^p (\mathbb{F})$ the space of $\mathbb{F}$-adapted  processes $Y$ such that

$$
\|Y\|^p_{\mathbf{B}^p}:=\mathbb{E}\left(\sup_{0\le t\le T}\|Y(t)\| \right)^p < \infty.
$$

\begin{definition}\label{continuouscontrolled}
A continuous \textbf{controlled functional} is a map $X:U_0\rightarrow L^p_a(\mathbb{P}\times Leb)$ for some $p\ge 1$, such that for each $t\ge 0$ and $u\in U_0$, $\{X(s,u); 0\le s\le t\}$ depends on the control $u$ only on $(0,t]$ and $X(\cdot,u)$ has continuous paths for each $u\in U_0$.
\end{definition}

 In many interesting examples, a continuous functional is in fact a Wiener functional. Our discretization technique will require that we consider ({\it payoff}) functionals also in the space of discontinuous trajectories. We shall denote by  $\mathbf{D}^n_T$ the space of c\`adl\`ag $\R^n$-value functions on the interval $[0,T]$ (recall that the French acronym c\`adl\`ag refers to continuity on the right with limit on the left). We endow the linear space $\mathbf{D}^n_T$ with the norm of uniform convergence on $[0,T]$.

We shall assume the following hypotheses which concern a sort of continuity of the {\it payoff} functional $\xi$ and of the controlled functional $X(\cdot, u)$
with respect to appropriate topologies.

\noindent \begin{assum}\label{AssumpA1} (H\"older continuity of $\xi$) There exists $\gamma \in (0,1]$ and a constant $C>0$ such that

\begin{equation*}\label{descA1}
|\xi(f) - \xi(g)|\le C \left(\sup_{0\le t\le T}\|f(t)-g(t)\|\right)^\gamma
\end{equation*}
for every $f,g\in \mathbf{D}^n_T$.

\end{assum}

\noindent \begin{assum}\label{AssumpB1} (Continuity of $u\mapsto X(\cdot, u)$) There exists a constant $C$ such that

\begin{equation*}\label{LipL2}
\|X(\cdot,u)-X(\cdot,v)\|^2_{\mathbf{B}^2(\mathbb{F})}\le C \mathbb{E}\int_0^T\|u(s)-v(s)\|^2_{\mathbb{R}^r}ds
\end{equation*}
for every $u,v \in U_0$.
\end{assum}

From now on, we are going to fix a controlled  functional $X:U^T_0\rightarrow\mathbf{B}^2(\mathbb{F})$.
For a given functional $\xi:\mathbf{D}^n_T\rightarrow\mathbb{R}$, we denote

\begin{equation*}\label{actionmap}
\xi_X(u): = \xi\big(X(\cdot,u)\big); u\in U_0
\end{equation*}
and the {\it value process}:
\begin{equation}\label{valuepdef}
V(t,u):=\text{ess}~\sup_{v\in U^T_t}\mathbb{E}\Big[\xi_X(u\otimes_t v)|\mathcal{F}_t\Big];0\le t< T,u \in U_0,
\end{equation}
where $V(T,u) := \xi_X(u)$~a.s and the process $V(\cdot, u)$ has to be viewed backwards. Throughout this paper, in order to keep notation simple, we omit the dependence of the value process in (\ref{valuepdef}) on the controlled functional $X$ and on the payoff $\xi$, in such a way that   $V$ means the mapping $V:U_0\rightarrow L^1(\mathbb{P}\times Leb)$.

Now,  we introduce what we call a \textit{controlled embedded discrete structure}, see, e.g.  \cite{LOS}. The heuristic is to view a controlled functional $u\mapsto Y(\cdot,u)$ as a family of simplified models which one has to build in order to extract the relevant information in order to get a concrete description of value processes and the construction of their associated $\epsilon$-optimal controls.

The discretization procedure will be based on a class of pure jump processes driven by suitable waiting times which describe the local behaviour of the Brownian motion. We briefly recall the basic properties of this skeleton (see among others, \cite{LOS} and references therein). We start by constructing a sequence of stopping times $\mathcal{T} := \{T^k_n; n\ge 0\}$  which is the basis for our discretization scheme. We set $T^{k}_0:=0$ and

\begin{equation*}\label{stopping_times}
T^{k}_n := \inf\{T^{k}_{n-1}< t <\infty;  \parallel  B(t)-B(T^{k}_{n-1}) \parallel_{\infty} = \epsilon_k\}, \quad n \ge 1,
\end{equation*}
where $\sum_{k\ge 1}\epsilon_k^2 < \infty$ and $\parallel\cdot \parallel_{\infty}$ is the $\ell_ {\infty}$ norm in $\R^n$. This implies that
\begin{equation*}\label{deltaMULT}
\Delta T^k_{n}:=T^k_n-T^k_{n-1}=\min_{j\in\{1,2,\dots, d\}}{\{\Delta T^{k,j}_n\}}~a.s
\end{equation*}
where, coordinate-wise, for $j=1, \ldots, d$

\begin{equation*}
\Delta T^{k,j}_n :=  \inf\{0< t <\infty;  | B^j(t+T^{k}_{n-1}) - B^j(T^{k}_{n-1}) | = \epsilon_k\}, \quad n \ge 1.
\end{equation*}
Then, we define the discretization of the Brownian motion $A^k :=(A^{k,1}, \cdots , A^{k,d})$ by

\begin{equation*}\label{rw}
A^{k,j} (t) := \sum_{n=1}^{\infty} \left(  B^j(T^{k}_{n}) -B^j(T^{k}_{n-1}    \right) \mathds{1}_{\{T^{k}_n\leq t \}};~t\ge0, ~ j=1, \ldots , d,
\end{equation*} for every $k\ge 1$.

The completion (including subsets of null sets) of the filtration generated by $A^{k}$ is denoted by $\mathbb{F}^{k} = (\mathcal{F}^{k}_t)_{t\ge 0}$. It is  naturally a sub-filtration of the Brownian motion's original filtration. Moreover:

$$
{\mathcal{F}}^k_t \cap\{T^k_n \le t < T^k_{n+1}\} = {\mathcal{F}}^k_{T^k_n}\cap \{T^k_n \le  t < T^k_{n+1}\}; t\ge 0
$$


\begin{definition}\label{discreteskeleton}
The structure $\mathscr{D} = \{\mathcal{T}, A^{k}; k\ge 1\}$ is called a \textbf{discrete-type skeleton} for the Brownian motion.
\end{definition}

By strong Markov property, we observe that (see, e.g. \cite{LEAO_OHASHI2018.1}):

\begin{enumerate}
  \item The jumps $\Delta A^{k,j}(T^k_n); n=1, 2, \ldots$ are independent and identically distributed (iid).
  \item The waiting times $\Delta T^k_n; n=1, 2,\ldots$ are iid random variables in $\mathbb{R}_+$.
  \item The families $(\Delta A^{k,j}(T^k_n); n=1, 2, \ldots)$ and $(\Delta T^k_n; n=1, 2,\ldots)$ are independent.
\end{enumerate}
Moreover, it is immediate that $A^{k,j}$ is a square-integrable $\mathbb{F}^k$-martingale for each $j=1,\ldots, d$.

This structure is the basic settings for many approaches of near-optimal problems. In particular, in the next section we are going to describe how we apply this framework in our dynamic programming technique.

\section{State of the art and applications}

We introduce the following  subclasses of locally constant control functions  $U^{k,T^k_n}_{T^k_\ell}\subset U_{T^k_\ell}^{T^k_n}$ with $ 0\le \ell < n < \infty$. Namely,  $U^{k,T^k_n}_{T^k_\ell}$ is the set of $\mathbb{F}^k$-predictable processes of the form

\begin{equation}\label{controlform0}
v^k(t) = \sum_{j=\ell+1}^{n}v^{k}_{j-1}\mathds{1}_{\{T^k_{j-1}< t\le T^k_j\}}; \quad T^k_\ell < t \le T^k_n,
\end{equation}
where for each $j=\ell+1, \ldots, n$, $v^k_{j-1}$ is an $A$-valued $\mathcal{F}^k_{T^k_{j-1}}$-measurable random variable. To keep notation simple, we use the shorthand notations

\begin{equation}\label{uknm}
U^{k,n}_\ell: = U^{k,T^k_n}_{T^k_\ell}; 0\le \ell < n,
\end{equation}
where $U^{k,\infty}_\ell$ is the set of all controls $v^k:~]]T^k_\ell, +\infty[[\rightarrow A$ of the form

\begin{equation*}\label{ukinf}
v^k(t) = \sum_{j\ge \ell+1}v^{k}_{j-1}\mathds{1}_{\{T^k_{j-1}< t\le T^k_j\}}; \quad T^k_\ell < t,
\end{equation*}
where $v^k_{j-1}$ is an $A$-valued $\mathcal{F}^k_{T^k_{j-1}}$-measurable random variable for every $j\ge \ell+1$ for an integer $\ell \ge 0$.

We also use a shorthand notation for $u^k\otimes_{T^k_\ell} v^k$ as follows: for $u^k\in U^{k,n}_0$ and $v^k\in U^{k,n}_{\ell}$ with $\ell < n$, we write

\begin{equation*}\label{Kconcatenation}
(u^k\otimes_\ell v^k): = (u^k_0, \ldots, u^k_{\ell-1}, v^k_\ell, \ldots, v^k_{n-1} ).
\end{equation*}
This notation is consistent because $u^k\otimes_{T^k_\ell} v^k$ only depends on the list of variables $(u^k_0, \ldots, u^k_{\ell-1}, v^k_\ell, \ldots, v^k_{n-1} )$ whenever $u^k:~]]0,T^k_n]]\rightarrow A$ and $v^k:~]]T^k_\ell,T^k_n]]\rightarrow  A$ are controls of the form (\ref{controlform0}) for $\ell < n$.


Let us now  introduce the analogous concept of controlled functional but based on the filtration $\mathbb{F}^k$. For this purpose, we need to introduce some further notations. Let us define

\begin{equation*}\label{ektdef}
e(k,t):= \Big\lceil \frac{\epsilon^{-2}_k t}{\chi_d}\Big\rceil; 0\le t\le T,
\end{equation*}
where $\lceil x\rceil$ is the smallest integer greater or equal to $x\ge 0$ and

\begin{equation*}\label{kappaDEF}
\chi_d:=\mathbb{E}\min\{\tau^1, \ldots, \tau^d\},
\end{equation*}
where $(\tau^j)$ is the stopping time $\inf\{t>0; |W^j (t)|=1\}$ for each $j=1, \ldots ,  d$.  We have the following crucial uniform convergence:
\begin{lemma} When $k$ goes to infinity, the random variable  $T^k_{e(k,s)}$ a.s. converges uniformly to $s$ in the interval $[0,t]$, i.e.  for any $0\leq t \leq T$
\[
\lim_{k\rightarrow \infty} \sup_{0\leq s \leq t} |T^k_{e(k,s)}-s| = 0 \ \ \ a.s.
\]
Moreover, the convergence $T^k_{e(k,t)}\rightarrow t$  holds also in $L^p  (\mathbb{P})$ for any $1\leq p \leq \infty$.

\end{lemma}
\proof: The uniform convergence holds using analogous arguments of \cite[Lemma 2.2]{Khoshnevisan and Lewis}, just adapting easily the argument there for $\epsilon_k = 2^{-n/2}$ to our general $\epsilon_k \in \ell_ 2$. Convergence in $L^p$ holds by Lebesgue convergence theorem.

\eop

Let $O_T(\mathbb{F}^k)$ be the set of all step-wise constant $\mathbb{F}^k$-adapted processes of the form

$$Z^k(t) = \sum_{n=0}^\infty Z^k(T^k_n)\mathds{1}_{\{T^k_n\le t\wedge T^k_{e(k,T)} < T^k_{n+1}\}}; 0\le t\le T,$$
where $Z^k(T^k_n)$ is $\mathcal{F}^k_{T^k_{n}}$-measurable for every $n\ge 0$ and $k\ge 1$.

Let us now present two concepts which will play a key role in the main result.
\begin{definition}\label{GASdef}
A \textbf{controlled embedded discrete structure} $\mathcal{Y} = \big((Y^k)_{k\ge 1},\mathscr{D}\big)$ consists of the following objects: a discrete-type skeleton $\mathscr{D}$ and a map $u^k\mapsto Y^{k}(\cdot,u^k)$ from $U^{k,e(k,T)}_0$ to $O_T(\mathbb{F}^k)$ such that

\begin{equation*}\label{antiprop}
Y^{k}(T^k_{n+1},u^k)~\text{depends on the control only at}~(u^k_0, \ldots, u^k_n),
\end{equation*}
for each integer $n\in \{0,\ldots,e(k,T)-1\}$.
\end{definition}

\begin{definition}\label{GASdefSTRONG}
A \textbf{strongly controlled functional} is a pair $(X,\mathcal{X})$, where $X$ is a controlled functional and $\mathcal{X} = \big((X^k)_{k\ge 1},\mathscr{D}\big)$ is a controlled embedded discrete structure such that the random variable $\big\{\|X^{k}(\phi)\|_\infty; \phi \in U^{k,e(k,T)}_0\big\}$ is uniformly integrable for each $k\ge 1$ and
\begin{equation*}\label{keyassepsilon}
\lim_{k\rightarrow +\infty}\sup_{\phi\in U^{k,e(k,T)}_0}\mathbb{E}\|X^{k}(\phi) - X(\phi)\|^\gamma_{\infty}=0,
\end{equation*}
for $0 < \gamma \le 1$.
\end{definition}

The concepts of controlled embedded discrete structures and strongly controlled  functionals are nonlinear versions of the structures analysed in \cite{LEAO_OHASHI2018.2}. The typical example we have in mind of a strongly controlled (Wiener) functional is a controlled state which drives a stochastic control problem. See example in the Section \ref{Numerical Example} below.

Throughout this section, we are going to fix a controlled embedded discrete structure

\begin{equation*}\label{controlledstate}
u^k\mapsto X^{k}(\cdot,u^k)
\end{equation*}
and we set

\begin{equation*}\label{Kactionmap}
\xi_{X^k}(u^k):=\xi\big(X^k(\cdot, u^k)\big),
\end{equation*}
for $u^k\in U^{k,e(k,T)}_0$. We assume $\big\{\|X^{k}(\phi)\|_\infty; \phi \in U^{k,e(k,T)}_0\big\}$ is uniformly integrable for each $k\ge 1$. We then define the discrete value process:

\begin{equation*}\label{discretevalueprocess}
V^{k}(T^k_n, u^k):=\esssup_{\phi^k\in U^{k,e(k,T)}_n}\mathbb{E}\Big[\xi_{X^k}(u^k\otimes_n\phi^k)\big|\mathcal{F}^k_{T^k_{n}}\Big]; \ \ \ \ \  n=1,\ldots, e(k,T)-1,
\end{equation*}
with boundary conditions

$$V^k(0):=V^k(0,u^k):=\sup_{\phi^k\in U^{k,e(k,T)}_0}\mathbb{E}\big[\xi_{X^k}(\phi^k)\big],\quad V^{k}(T^k_{e(k,T)}, u^k): = \xi_{X^k}(u^k).$$

This naturally defines the map $V^k:U^{k,e(k,T)}_0\rightarrow O_T(\mathbb{F}^k)$ with jumps $V^k(T^k_n,u^k); n=1, \ldots, e(k,T)$ for $u^k\in U^{k,e(k,T)}_0$. One should notice that $V^k(T^k_n, u^k)$ only depends on $u^{k,n-1}:=(u^k_0, \ldots, u^k_{n-1})$ so it is natural to write

$$V^k(T^k_n, u^{k,n-1}) :=V^k(T^k_n, u^k); u^k\in U^{k,e(k,T)}_0, 0\le n\le e(k,T),$$
with the convention that $u^{k,-1}:=\mathbf{0}$. By construction, $V^k$ satisfies the definition of a controlled embedded discrete structure.

\

We present the main theorem in this abstract setting, which naturally includes the original optimization problem of the original model (\ref{eq: model}). The discretization structure constructed so far plays a fundamental role:

\

\begin{theorem}\label{Thm: PD}
Let $\mathcal{X}$ and $\xi:\mathbf{D}^{n}_T\rightarrow\mathbb{R}$ be respectively a strongly controlled functional and a (payoff) functional which satisfies the continuity Assumptions (1) and (2) of the last Section.
Then:

\begin{description}
\item[1)] The value function $(V^{k}_j)_{j=0}^{e(k,T)}$ associated with $X^k$ is the unique solution of

\begin{equation} \label{keyrecursion}
\esssup_{u^k \in U^{e(k,T)}_j}\mathbb{E}\left(V^k(T^k_{j+1},\theta\otimes_j u^k)\big|\mathcal{F}^k_{T^k_j}\right)-V^k(T^k_j,\theta)  =  0, ~ a.s. \quad j=e(k,T)-1, \ldots, 0,
\end{equation}
for all $\theta \in U_0^{j}$,  with boundary condition
\[
V^k (T^k_{e(k,T)},u) =\xi_X^k\left( u\right) ~~~~\mbox{for } ~~~ u \in U^{e(k,T)}_0.
\]

\item[2)] We have
\begin{equation*}
    \lim_{k\rightarrow \infty}\Big|V^k(0)-\sup_{u\in U_0^T}\mathbb{E}\left(\xi_{X}(u)\right)\Big|=0
\end{equation*}

\item[3)] There exist control functions $u^{k, \star}_j:\Omega\times [T^k_{j},T^k_{j+1}]\rightarrow\mathbb{R}$ such that

\begin{equation*}
 V^{k}(T^k_j,\theta) = \mathbb{E}\left(V^k(T^k_{j+1},\theta\otimes_j u^{k, \star}_j)\Big|\mathcal{F}^k_{T^k_n}\right)
\end{equation*}
where $j=e(k,T)-1,\ldots, 1$.

\item[4)] Given $\epsilon>0$,
\begin{equation*}\label{keyorig}
\mathbb{E}\left(\xi_{X^k}(u^{k,\star})\right)\geq \sup_{u\in U_0^T}\mathbb{E}\left(\xi_{X}(u)\right)-\epsilon
\end{equation*}
for every $k$ sufficiently large, with $u^{k,\star}=(u^{k,\star}_0,\dots,u^{k,\star}_{e(k,T)-1})$.
\end{description}

\end{theorem}

Note that in the Theorem above, we have transformed the infinite-dimensional optimization problem in a finite-dimensional algorithm which can be numerically implemented. The near-optimal problem of equation (\ref{cp1}) is then solved for a reasonable wide range of systems.

\noindent {\bf Comments on the proof and numerical algorithms:} Initially note that formula (\ref{keyrecursion}) is just the algorithm of the dynamic programming. For its implementation one has to  know explicitly the model of the controlled system, i.e. the parameters of equation (\ref{eq: model}). Nevertheless, numerically speaking, mind that if the  model is unknown (e.g. typical situation in finance), one still can use machine learning techniques to obtain an optimal control and, consequently, the value function based on data and on the simulation of the environment.

For simplicity, for any control function $u^k_ j$ define
$$Q^{k}(T^k_j,\theta) = \mathbb{E}\left(V^k(T^k_{j+1},\theta\otimes_j u^{k}_j)\Big|\mathcal{F}^k_{T^k_n}\right).$$

We have several approaches on solving the dynamic programming equation. Say, given the model, one can try to solve it directly and obtain a closed solution for optimal control. For example, if the controlled system is linear, the instantaneous cost function and the terminal cost are quadratic functions, it is possible to determine the solution of the dynamic programming equation using the Riccati equation. In general, the numerical algorithms to solve the dynamic programming equation are based on:

\begin{itemize}
  \item Choose the discretization level $k$ using the convergence rate (see details  in \cite{LOS}).
  \item Approximate $Q^k_j $ function into each step $j$ by calculating the conditional expectation. For example, use regression methods with Monte Carlo simulation and machine learning techniques.
  \item Once we get an approximation $\hat{Q}^k_j $ to the function $Q^k_j $, the optimal control $\hat{u} (\mathbf{o}^k_{j}) $ obtained by exhaustive search if $ A$ is discrete, or by algorithms based on the gradient method via e.g. neural networks, when $A$ is a compact with open interior.
\end{itemize}

\eop

For a full proof of the theorem, see \cite{LOS}. Recently, several numerical methods based on deep learning techniques have been proposed to solve the PD equation, see \cite{HW,HPBL,HPBLN}. In this context, Huré et al. \cite {HPBL} propose machine learning based methods to solve the dynamic programming equation. In general, under the assumption that the controlled system is Markovian, these authors propose an approach for optimal control through deep learning techniques and calculate the value function via Monte Carlo simulation. Theorem \ref{Thm: PD} above proposes a generalization of these algorithms for the non-Markovian case.

\subsection{Numerical Example} \label{Numerical Example}

Let us now present a simple example to illustrate the theory presented so far. We intent to show the power of the numerical technique in a toy model in such a way that one can compare the results of the simulation with the unique theoretical solution. For this purpose, we choose the example of hedging in a Black-Scholes model:   For given real parameters $c$ and $K>0$, consider the cost functional

$$\varrho_c(x,y,K):=(c+ x -\max(y-K,0))^2; \ (x,y)\in \mathbb{R}^2.$$
Assume that we have the classical geometrical equation:
\begin{equation*}
\begin{split}
dS(t)&=S(t)\left(\mu dt+\sigma dB(t)\right)\\
\end{split}
\end{equation*}
where, for simplicity, we assume drift $\mu=0$, final time $T=1$. Then, the problem is to find, among all possible control functions $\phi\in U^T_0$

\begin{equation*}\label{numexfinal}
\inf \mathbb{E}\Big[\varrho_c\big(X(T,\phi),S(T), K\big)\Big]\quad ,
\end{equation*}
where the controlled state space is given by:

\begin{equation*}  \label{Eq: geometrical}
X(t,\phi)=\int_0^t \phi(r)dS(r); \ \ \ \phi\in U^T_0, 0\le t\le T.
\end{equation*}
The control functions represent the absolute percentages of the securities $S$ which an investor holds at time $t\in [0,T]$ when, in our notation $a=1$. It is well-known that there exists a unique choice of $(c^*,\phi^*) \in \mathbb{R}\times U^T_0$ such that

\begin{equation*}
\inf_{(c,\phi)\in \mathbb{R}\times U^T_0}\mathbb{E}\Big[\varrho_c\big(X(T,\phi),S(T), K\big)\Big] = \mathbb{E}\Big[\varrho_{c^*}\big(X(T,\phi^*),S(T), K\big)\Big]=0,
\end{equation*}
whose solution is given by

$$c^*=S_0\Phi(d_1)-\frac{K}{B_T}\Phi(d_2),$$
where

$$ d_1=\frac{\log\left(\frac{S(0)}{K}\right)+\frac{\sigma^2}{2}T}{\sigma\sqrt{T}},\quad d_2=d_1-\sigma\sqrt{T},$$
and $\Phi$ is the cumulative distribution function of the standard Gaussian variable. We recall $c^*$ is the price of the option and $\phi^*$ is the so-called delta hedging which can be computed by means of the classical PDE Black-Scholes as a function of $\Phi$. We set $S^1(0)=49, \sigma=0.2$, $K=55$ and $\epsilon_k = 2^{-k}$. The discretization is given by $$X^{k}(t,\phi^k) = \mathbb{X}^k(t\wedge T^k_{e(k,T)},\phi^k),\quad \mathbb{X}^{k}(t,\phi^k)= \int_0^{\bar{t}_k}\phi^k(r)dS^{k}(r); \phi^k\in U^{k,e(k,T)}_0,$$
and $S^{k}(t)  = \mathbb{S}^{k}(t\wedge T^k_{e(k,T)})$. To shorten notation write $S^{k}_n = S^{k}(T^k_n);  0\le n \le e(k,T)$ and  $\mathbb{E}_{k,n}$ for a conditional expectation with respect to $\mathcal{F}^k_{T^k_n}$. At first, for a given $c\in \mathbb{R}$, we apply the algorithm described above to get a Monte Carlo optimal control approximation $\phi^{*,k} = (v^{*,k}_0, \ldots, v^{*,k}_{m-1})$. In this particular case, we can analytically solve  and the optimal control is given by: $v^{\star,k}_{m}=0$ and for $1\le i\le m$, we have

\begin{equation*}
\begin{split}
v^{\star,k}_{m-i}&=-\frac{\mathbb{E}_{k,m-i}\left[\sum_{j=0}^{i-1}\left(v^{\star,k}_{m-j}\Delta S^{k}_{m+1-j} -H\right)\Delta S^{k,1}_{m-i+1} \right]}{\left(S^{k,1}_{m-i}\right)^2\sigma_1^2\epsilon^2_k}\\
\end{split}
\end{equation*}
with $H=\max{\{(S(T)-K),0\}}$. The estimated value $c^{k,*}$ is computed according to

\begin{equation}\label{cstar}
\begin{split}
c^{k,*}\in \argmin_{c\in \mathbb{R}}\mathbb{E}\left[\varrho_c\left(X^k(T,\phi^{k,*}),S^{k}(T\wedge T^k_{e(k,T)}), K\right)\right].
\end{split}
\end{equation}

In other words, $\mathbb{E}\left[\varrho_{c^{k,*}}\left(X^k(T,\phi^{k,*}),S^{k,1}(T\wedge T^k_{e(k,T)}), K\right)\right]=0$, since the minimum of the function (\ref{cstar}) is zero. Table \ref{001tab} presents a comparison between the true call option price $c^*$ and the associated Monte Carlo price $c^{k,*}$. Figure 1 presents the Monte Carlo experiments for $c^{k,*}$ with $k=3$. The number of Monte Carlo iterations in the experiment is $2 \times 10^4$.

\begin{figure}[!h]
  \centering
  \includegraphics[width=0.7\textwidth]{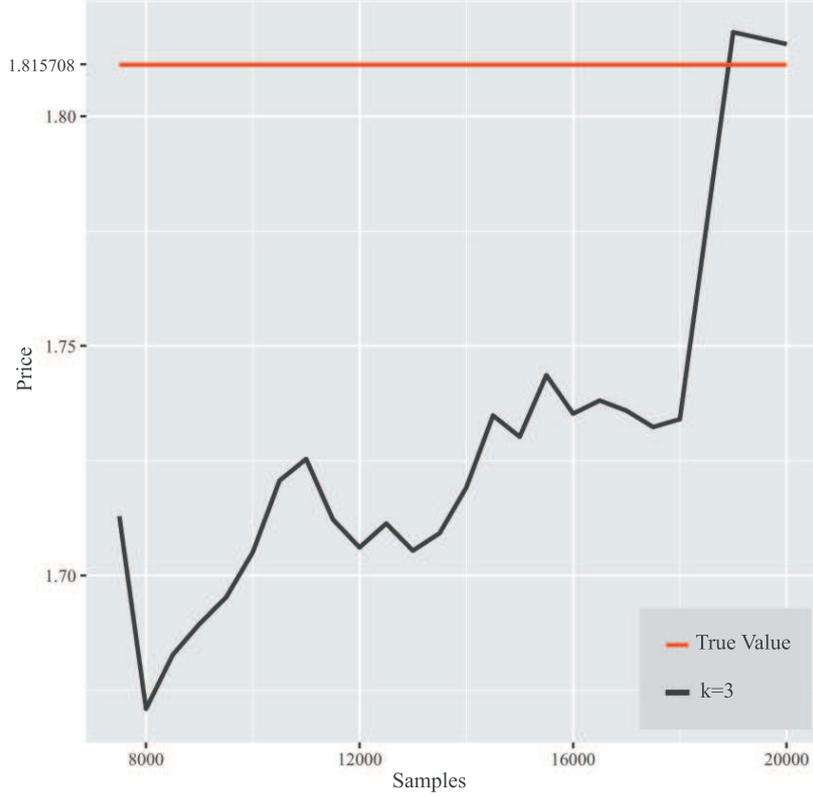}
  \caption{Monte Carlo experiments for $c^{k,*}$ with $\epsilon_k = 2^{-k}$}\label{1fig}
\end{figure}

\begin{table}[!h]
\caption{Comparison between $c^{*}$ and $c^{k,*}$ for $\epsilon_k = 2^{-k}$}
\begin{tabular}{lccccc}
  \hline
  \textbf{Result} & \textbf{Mean Square Error} & \textbf{True Value} & \textbf{Difference} & \textbf{\% Error}\\
  \hline
    1.815708  & 0.008660 & 1.811209 & 0.004499 & 0.002477\%\\
  \hline
\end{tabular}
\label{001tab}
\end{table}

\section{Open problems}

In this section we describe open problems which are natural extended situations where the near-optimal stochastic control problem described in the previous sections are relevant. One of them is the more geometrical case of Lie groups and their corresponding quotients (homogeneous spaces) which include spheres, torus, hyperboloid models and many others that are the natural framework for many interesting applications.  Another extended direction to be explored is the case of noise with jumps, for example semimartingales with jumps in general, like \cite{KPP} or, in particular, Lévy noise as e.g. \cite{D.Applebaum}. Here we intend to describe these open problems and some of the difficulties on dealing with them.

\subsection{Lie groups}
Besides the pure geometrical motivation on extending the control settings to Lie groups, we point out some other applied importance: 1) this theory encapsulates the multiplicative approach, in the sense that the products here can be considered as the usual product of square matrices; 2) Lie groups are the appropriate frame to work with linear systems: in many of these cases $G$ is the group (or a closed subgroup ) of positive determinant $n\times n$-matrices $Gl^+(n, \R )$; 3)  this approach also includes the framework for any homogeneous space via quotient by closed subgroups, e.g. $n$-dimensional spheres, torus, Grasmannian, projective spaces, hyperboloid model and many others. In all these cases, the properties of the Lie group theory are important tools in the analysis and interpretation of the dynamics:  decompositions (Iwasawa, polar, eingenvectors etc), invariance by translations, adjoints, geometrical structures of fibre bundles, connections, local diffeomorphism with to the corresponding Lie algebra, and many more. A vast and wide-ranging literature are available from applied to more theoretical approach. Just to mention few of them, see e.g. from the classical \cite{Chevalley}, the well known \cite{Warner}, the more introductory \cite{Baker} or the more recent \cite{San Martin}, including all references therein.

Let $G$ be a connected Lie group with  its corresponding Lie algebra  $\mathfrak{g}$, identified with the tangent space $T_e G$, where $e$ is the identity of $G$. The dynamics of the controlled trajectories $X^u(t)\in G$, with initial condition $X^u(0)=e$ in this context is described by the right invariant vector fields:

\begin{equation} \label{Eq: control in Lie group}
dX^u (t) = d(R_{X^u (t)})_e\  \alpha(t, X^u_ t, u(t))\, dt + d (R_{X^u (t)})_e\  \sigma(t, X^u_ t, u(t))\, dZ(t),
\end{equation}
where $\alpha $ and $\beta$ are $\mathfrak{g}$-valued functions, $d(R_{X^u (t)})_e: T_eG \rightarrow T_{X^u (t)}G$ is the derivative at the identity of the right translation $R_{X^u} (t):  g\mapsto g  X^u (t)$. As before, the functions $\alpha $ and $\beta$ depend on time $t$, the past trajectory $s\mapsto X^u (s)$ with $s\in [0,t]$ and the control function (bounded measurable) $u$. If the Lie group $G$ is a subgroup of matrices, equation (\ref{Eq: control in Lie group}) can be written with a much more familiar notation
\begin{equation*} \label{Eq: control matrices}
dX^u (t) =   \alpha(t, X^u_ t, u(t)) \cdot  X^u (t) \, dt +   \sigma(t, X^u_ t, u(t))\cdot  X^u (t)\  dZ(t),
\end{equation*}
where $\cdot$ stands for the usual product of square matrices.

 In this new context, the cost (or payoff) functional are defined on $C_T$, the set of continuous trajectories in $G$. Hence, as before, for a certain final time $T>0$ one is looking for the optimal performance given by
\begin{equation*}
\sup_{u\in U^T} \E \left[  \xi (X^u)\right].
\end{equation*}
As before, this optimal control function may not exist and the problem turns into a near-optimal control function. Previous Assumptions (1) and (2) can be consider naturally in this context. The discrete type skeleton is again a helpful tool here: note that the Lie group structure allows a special kind of discretization of the dynamics given by a product of random matrices $M_k\cdot  \ldots M_2 \cdot M_1$, where each $M_j$, $j=0,1, \ldots, k$ is an exponential in the Lie group $G$ with $M_0=e$ a.s.. Hence the main point here is that controlling $X^u (t)$ is somehow the same as controlling the distributions of $M_j's$. To understand fully this multiplicative approach, i.e. in each sense its structure can provide new properties to the solution demands good answers to some questions like: how to translate the cost functional into this product of matrices? What are the distributions of these random matrices? In what sense controlling $X^u (t)$ is the same as controlling the distributions of $M_j's$? One needs a comparison of the dynamical programming in the additive case, as stated in Theorem \ref{Thm: PD}, with a DP algorithm here which will start from left to the  right in the product of the $M_j 's$. For many applications it is also relevant to understand how this approach behaves with respect to quotients, i.e. for control systems in homogeneous space.

\subsection{Discontinuous noise and geodesic jumps}

 Another interesting open problem related to optimal stochastic control relates to discontinuous noise. In particular, we shall consider semimartingales with jumps such that the  trajectories are càdlàg. This is the classical case, see e.g. Kurtz, Pardoux and Prother \cite{KPP}, which includes, for example the Lévy noise, Applebaum \cite{D.Applebaum}, Protter \cite{Protter}, Oksendal and Sulem \cite{Oksendal}, among others.  In this context $Z = \{ Z_t, t\geq 0\}$ is a $k$-dimensional semimartingale, $[Z,Z]$ denotes the matrix of quadratic variation of $Z$ which can be decomposed into $[Z,Z] = [Z,Z]^c + [Z,Z]^d$, where $[Z,Z]^c$ and $[Z,Z]^d$ represent the continuous and purely discontinuous parts respectively. This kind of noise represents the appropriate model for many applications, in fact, just note that the discretization in time of any system driven by Brownian noise (with or without Poisson perturbation) turns out to be of this kind. The open questions here relates to the behaviour of the DP algorithm of Theorem \ref{Thm: PD} with respect to jumps in the noise.

The canonical Marcus equation is an interesting approach for the dynamics of jump diffusions on manifolds. Heuristically, if $\Delta_s Z = Z(s) - Z(s^-)$ is the jump of the noise at time $s$, the corresponding jump in the dynamics is performed by a ficticious time-one jump along the deterministic flow of the corresponding vector field multiplied by $\Delta_s Z$. For more details we refer to \cite{KPP} where many other properties are also presented. The second geometrical aspect that we propose in this open problem is to perform the jump along the ficticious time-one to be performed along the geodesic starting at the point $X(s^-)$ in the direction of vector field at this point multiplied by $\Delta_s Z$. We are going to make it clear in the next paragraph.

Let $M$ be an $m$-dimensional complete Riemannian manifold. Let $\alpha(t, X^u_t, u(t))$ and $\sigma(t, X^u_t, u(t))$ be smooth vector fields with bounded derivatives in $M$. We assume that the controlled process $X^u(t)\in M$ is modelled by the following geodesic Marcus equation:
\begin{equation}
dX^u(t) = \alpha(t, X^u_t, u(t))\ dt + \sigma(t, X^u, u(t)) \ \diamond dZ(t),
\label{eq1}
\end{equation}
which is defined in the following way:
$X^u(t)$ is a solution of equation (\ref{eq1}) if and only if for any test function $f\in C^2(M)$ we have:
\begin{eqnarray}
f(X^u(t)) - f(X^u(0)) &=& \int_0^t \alpha(s, X^u_s, u)s))ds + \int_0^t \sigma f(s, X^u_s, u(s))dZ(t) \nonumber \\
&+& \frac{1}{2}\int_0^t \sigma' \sigma(s, X^u_s, u(s)) f(X^u_s) d[Z,Z]_s^c \nonumber \\
&+& \sum_{0<s\leq t}\left\lbrace f(\gamma(X^u_{s-}, \sigma(s, X^u_s, u(s)) \Delta_s Z)) \right. \nonumber \\
&-& \left. f(X^u_{s-}) - \sigma(s, X^u_{s-}, u(s)) f(X^u_{s-})\Delta_s Z \right\rbrace.
\label{eq2}
\end{eqnarray}
Here $\gamma(X^u_{s-}, \sigma(s, X^u_s, u(s)) \Delta_s Z)$ is the unique geodesic starting at $X^u_{s-}$ in the direction of $\sigma(s, X^u_s, u(s)) \Delta_s Z \in T_{X^u_{s^-}}M$.
We remark some aspects. Firstly, note that one can easily includes a higher dimensional noise acting in more than one vector field in the sense that

\begin{eqnarray*}
\int_0^t \sigma f(s, X^u_s, u(s)) dZ(t) := \int_0^t \sigma^i f(s, X^u_s, u(s))dZ^i(t).
\end{eqnarray*}
Equally
\begin{eqnarray*}
\int_0^t \sigma' \sigma(s, X^u_s, u(s)) f(X^u_s) d[Z,Z]^c_s := \int_0^t \frac{D\sigma}{\partial x^i} \sigma^i(s, X^u_s, u(s)) f(X^u_s) d[Z^i,Z^j]^c_s,
\end{eqnarray*}
is a Stieltjes integral with respect to the continuous bounded variation processes $[Z^i,Z^j]^c_s$.

Secondly, note that if the test function $f$ is the identity in an Euclidean space and the rule of jumps $\gamma(X^u_{s-}, \sigma(s, X^u_s, u(s)) \Delta_s Z)$ is converted to the time-one flow of  $\sigma$, then equation (\ref{eq2}) turns back to be the classical Marcus equation as in \cite{KPP}. The proposed geodesic jumps corresponds to the natural generalization of the Itô equation with càdlàg noise in an Euclidian space to a Riemannian manifold. As expected, in this case we have a coalescing stochastic solution flow.

Finally, note that the subtracting terms on the summation in equation (\ref{eq2}) are already accounted in the It\^{o} integral on the right hand side. Applying Taylor's formula in $f(X^u_{s-}, \sigma(s, X^u_{s-}, u(s)) \Delta_s Z)$ and using the finite quadratic variation one sees that the series converges absolutely.

With this geodesic description of the dynamics, many properties has to be proved: say, existence, uniqueness, flows (of coalescing smooth maps),  Itô formula, Itô-Ventsel kind of formulas, probabilistic properties on the process, a kind of support theorem ($M$-invariance when $M$ is embededd in an Euclidean space with tangent vector fields) among many others.  In the context of the control problem, the challenge  is to answer the following question: how does the DP algorithm behaves in this context? What would be geometrical properties of the near-optimal control function $u^\star$ which satisfies $\displaystyle{\E \left[  \xi (X^{u^*})\right] > \sup_{u\in U^T} \E \left[  \xi (X^u)\right]} -\epsilon$, for any positive $\epsilon$ ? Dependence on geometrical parameters, curvature or even on the topology of $M$,  etc, looks quite promising and interesting issues.

\bibliographystyle{plain}

\end{document}